\documentclass[4pt, oneside]{article}

\usepackage[OT2, T1]{fontenc}
\usepackage[russian, english]{babel}
\usepackage{setspace}

\usepackage{enumerate}
\usepackage{amsmath,amssymb,amsthm,units,stmaryrd}
\usepackage{qtree,bussproofs}

\usepackage{yfonts}
\usepackage{units,stmaryrd}
\usepackage[colorlinks]{hyperref}
\usepackage[usenames,dvipsnames]{color}
\hypersetup{
  colorlinks,
  citecolor=Violet,
  linkcolor=Red,
  urlcolor=Blue}
\usepackage{graphicx}
\usepackage[all]{xy}
\newtheorem{theorem}{Theorem}[section]

\newtheorem{definition}[theorem]{Definition}
\newtheorem{lemma}[theorem]{Lemma}
\newtheorem{remark}[theorem]{Remark}
\newtheorem{corollary}[theorem]{Corollary}
\newtheorem{proposition}[theorem]{Proposition}

\usepackage{xspace}

\newcommand{\la}{\langle \,}
\newcommand{\ra}{\, \rangle}

\def\le{{\ell}}

\def\fmodels{\xymatrix{
\ar@{|=}[r]^{<\omega}&
}
}
\def\nmodels{\xymatrix{
\ar@{|=}[r]^{N}&
}
}
\def\<{\left <}

\def\>{\right >}

\DeclareSymbolFont{AMSb}{U}{msb}{m}{n}
\DeclareMathSymbol{\N}{\mathbin}{AMSb}{"4E}
\DeclareMathSymbol{\Z}{\mathbin}{AMSb}{"5A}
\DeclareMathSymbol{\R}{\mathbin}{AMSb}{"52}
\DeclareMathSymbol{\Q}{\mathbin}{AMSb}{"51}
\DeclareMathSymbol{\I}{\mathbin}{AMSb}{"49}

\newcommand{\bt}{\begin{theorem}}
\newcommand{\et}{\end{theorem}}
\newcommand{\bl}{\begin{lemma}}
\newcommand{\el}{\end{lemma}}

\newcommand\remove[1]{}

\newcommand{\On}{{\sf On}}
\newcommand{\TPr}{\textbf{\foreignlanguage{russian}{TS}C}}
\newcommand{\al}{\alpha}
\newcommand{\be}{\beta}
\newcommand{\FLE}{\mathbb{F}}
\newcommand{\MNF}[2]{\la n_{#1}^{\al_{#1}} \ra \top \land \ldots \land \la n_{#2}^{\al_{#2}} \ra \top }
\newcommand{\suc}{\text{Succ}}
\newcommand{\Lim}{\text{Lim}}

\def\le{{\ell}}

\def\I{{\ensuremath{\mathcal{J}}}\xspace}

\begin{document}

  \title{Relational Semantics for the Turing Schmerl Calculus}
  \author{Eduardo Hermo Reyes \footnote{\href{mailto:ehermo.reyes@ub.edu}{ehermo.reyes@ub.edu}} \\ Joost J. Joosten \footnote{\href{mailto:ehermo.reyes@ub.edu}{jjoosten@ub.edu} } \\ {\small University of Barcelona}\\ {\small Department of Philosophy}  }
\maketitle

  \begin{abstract}
In \cite{HermoJoosten:2017:TSCArith} the authors introduced the propositional modal logic $\TPr$ (which stands for Turing Schmerl Calculus) which adequately describes the provable interrelations between different kinds of Turing progressions. The current paper defines a model $\mathcal{J}$ which is proven to be a universal model for $\TPr$. The model $\mathcal{J}$  is a slight modification of the intensively studied $\mathcal{I}$ : Ignatiev's universal model for the closed fragment of G\"odel L\"ob's polymodal provability logic $\textbf{GLP}$.
  \end{abstract}

\section{Introduction}

Turing progressions arise by iteratedly adding consistency statements to a base theory. Different notions of consistency give rise to different Turing progressions. In \cite{HermoJoosten:2017:TSCArith}, the authors introduced the system $\TPr$ that generates exactly all relations that hold between these different Turing progressions given a particular set of natural consistency notions. The system was proven to be arithmetically sound and complete for a natural interpretation, named the \emph{Formalized Turing progressions} (FTP) interpretation.  

In this paper we discuss relational semantics of $\TPr$ by considering a small modification on Ignatiev's frame, which is a universal frame for the variable-free fragment of Japaridze's provability logic $\textbf{GLP}$.

\section{Strictly positive signature}

$\TPr$ is built-up from a positive propositional modal signature  using \emph{ordinal modalities}. Let $\Lambda$ be a fixed recursive ordinal throughout the paper with some properties as specified in Remark \ref{remark:WhyEpsilonNumbers}. By ordinal modalities we denote modalities of the form $\la n^\al  \ra$ where $\al \in \Lambda$ for some fixed ordinal $\Lambda$ and $n \in \omega$ (named \emph{exponent} and \emph{base}, respectively). The set of formulas in this language is defined as follows:

\begin{definition}
By $\FLE$ we denote the smallest set such that:

\begin{enumerate}[i)]
\itemsep0em
\item $\top \in \FLE$;
\item If $\varphi, \, \psi \in \FLE\Rightarrow (\varphi \wedge \psi) \in  \FLE$;
\item if $\varphi \in \FLE, \ n < \omega$ and $\al < \Lambda \Rightarrow \la n^\al \ra \varphi \in \FLE$.
\end{enumerate}
\end{definition}

For any formula $\psi$ in this signature, we define the set of base elements occurring in $\psi$. That is:

\begin{definition}
The set of base elements occurring in any modality of a formula $\psi \in \FLE$ is denoted by ${\sf N\text{-}mod} (\psi)$. We recursively define ${\sf N\text{-}mod}$ as follows:
\begin{enumerate}[i)]
\item ${\sf N\text{-}mod}(\top) = \emptyset$;
\item ${\sf N\text{-}mod}(\varphi \wedge \psi) = {\sf N\text{-}mod}(\varphi) \cup {\sf N\text{-}mod}(\psi)$;
\item ${\sf N\text{-}mod}(\la n^\al \ra \psi) = \{ n \} \cup {\sf N\text{-}mod}(\psi)$.	
\end{enumerate}
\end{definition}

\section{The logic $\TPr$}

In this section we introduce the logic $\TPr$ whose main goal is to express valid relations that hold between the corresponding Turing progressions. For this purpose we shall consider a kind of special formulas named \emph{monomial normal forms} which are used in the axiomatization of the calculus $\TPr$.\\

Monomial normal forms are conjunctions of monomials with an additional condition on the occurring exponents. In order to formulate this condition we first need to define the \emph{hyper-exponential} as studied in \cite{FernandezJoosten:2012:Hyperations}.

\begin{definition}
For every $n \in \omega$ the \emph{hyper-exponential} functions $e^n: \text{On} \rightarrow \text{On}$ are recursively defined as follows: $e^0 $ is the identity function, $e^1: \al \mapsto -1+\omega^\al$ and $e^{n+m}= e^n \circ e^m$.
\end{definition}

We will use $e$ to denote $e^1$. Note that for $\alpha$ not equal to zero we have that $e(\alpha)$ coincides with the regular ordinal exponentiation with base $\omega$; that is, $\alpha \mapsto \omega^\alpha$. However, it turns out that hyper-exponentials have the nicer algebraic properties in the context of provability logics.

\begin{definition}
The set of formulas in \emph{monomial normal form}, {\sf MNF}, is inductively defined as follows:
\begin{enumerate}[i)]
\item $\top \in {\sf MNF}$;
\item $\la n^\al \ra \top \in {\sf M}$, for any $n \, {<} \, \omega$ and $\al \, {<} \, \Lambda$;
\item \begin{tabbing}
	  if \hspace{0.2cm} \= $a) \ $ \= $\la n_0^{\al_0} \ra \top \wedge \ldots \wedge \la n_k^{\al_k} \ra \top \in {\sf MNF}$; \\
 	  \	\\	     
	     \> $b) \ $ \> \ $n < n_0$; \\
	  \ \\
	     \> $c) \ $ \> \ $\al$ of the form $e^{n_0 - n}(\al_0) \cdot (2 + \delta)$ for some $\delta < \Lambda$,
      \end{tabbing}
      \vspace{0.25cm}
then $\la n^\al \ra \top \wedge \la n_0^{\al_0} \ra \top \wedge \ldots \wedge \la n_k^{\al_k} \rangle \top \in {\sf MNF}.$
\end{enumerate}
\end{definition}

The derivable objects of $\TPr$ are \emph{sequents} i.e. expressions of the form $\varphi \vdash \psi$ where $\varphi, \, \psi \in \mathbb{F}$. We will use the following notation: by $\varphi \equiv \psi$ we will denote that both $\varphi \vdash \psi$ and $\psi \vdash \varphi$ are derivable. Also, by convention we take that for any $n$, $\la n^0 \ra \varphi$ is just $\varphi$.

\begin{definition} 
$\TPr$ is given by the following set of axioms and rules:\\

Axioms:
\begin{enumerate}

\item $\varphi \vdash \varphi, \ \ \ \varphi \vdash \top$;	

\item $\varphi \wedge \psi \vdash \varphi,  \ \ \ \varphi \wedge \psi \vdash \psi$;	\label{conjel}

\item Monotonicity axioms: $\la n^\al \ra \varphi \vdash \la n^\be \ra \varphi$, \ \ \ for $\be \, {<} \, \al$;	\label{mon}
		
\item Co-additivity axioms: $\la n^{\be + \al} \ra \varphi \equiv \la n^\al \ra \la n^\be \ra \varphi$;	 \label{coadditive}

\item Reduction axioms: $\la (n + m)^\al \ra \varphi \vdash \la n^{e^m (\al)} \ra \varphi$;	\label{omega}

\item Schmerl axioms: 
\[
\la n^\al \ra \big( \, \la n_0^{\al_0} \ra \top \ \wedge \ \psi \, \big) \equiv \la n^{e^{n_0 - n} (\al_0) \cdot (1 + \al)}  \ra \top \land \la n_0^{\al_0} \ra \top \ \land \ \psi
\] 
for $n\, {<} \, n_0$ and $\la n_0^{\al_0} \ra \top \ \wedge \ \psi \in {\sf MNF}$. \label{MS1}

\end{enumerate}
\

Rules:
\begin{enumerate}
\item If $\varphi \vdash {\psi}$ and $\varphi\vdash{\chi}$, then $\varphi\vdash {\psi}\wedge{\chi}$;	\label{r:1}
\item If $\varphi\vdash {\psi}$ and ${\psi}\vdash{\chi}$, then $\varphi\vdash{\chi}$;	\label{r:2}
\item If $\varphi\vdash {\psi}$, then $\la n^\al \ra \varphi\vdash \la n^\al \ra {\psi}$ ;	\label{r:3}
\item If $\varphi \vdash \psi$, then $\la n^\al \ra \varphi \, \land \, \la m^{\be + 1} \ra \psi \, \vdash \, \la n^\al \ra \big( \, \varphi \, \land \, \la m^{\be+ 1} \ra \psi \, \big)$ \ for $n \, {>} \, m$. \label{r:4}
\end{enumerate}
\end{definition}

It is worth mentioning the special character of Axioms (\ref{omega}) and (\ref{MS1}) since both axioms are modal formulations of principles related to Schmerl's fine structure theorem, also known as \emph{Schmerl's formulas} (see \cite{Schmerl:1978:FineStructure} and \cite{Beklemishev:2003:AnalysisIteratedReflection}).

\begin{remark}\label{remark:WhyEpsilonNumbers}
As we see in the axioms of our logic, they only make sense if the ordinals occuring in them are available. Recall that $\Lambda$ is fixed to be a recursive ordinal all through the paper. Moreover, some usable closure conditions on $\Lambda$ naturally suggest themselves. Since it suffices to require that for $n<\omega$ that $\alpha,\beta < \Lambda \ \Rightarrow \ \alpha + e^n(\beta) < \Lambda$, we shall for the remainder assume that $\Lambda$ is an $\varepsilon$-number, that is, a positive fixpoint of $e$ whence $e(\Lambda)=\Lambda=\omega^\Lambda$.
\end{remark}

In \cite{HermoJoosten:2017:TSCArith}, the authors proved that for any formula $\varphi$, there is a unique equivalent $\psi$ in monomial normal form. 

\begin{theorem} \label{MNFT}
For every formula $\varphi$ there is a unique $\psi \, {\in} \, {\sf MNF}$ such that $\varphi \equiv \psi$.
\end{theorem}

In virtue of the Reduction axioms, a formula $\psi \in {\sf MNF}$ may bear implicit information on monomials $\la n^\al \ra \top$ for $n \not \in {\sf N\text{-}mod}(\psi)$. The next definition is made to retrieve this information.

\begin{definition} \label{projection}
Let $\psi := \MNF{0}{k} \in {\sf MNF}$. By $\pi_{n_i} (\psi)$ we denote the corresponding exponent $\al_i$. Moreover, for $m \not \in {\sf N\text{-}mod}(\psi)$, with $n_k > m$, $\pi_m ( \psi)$ is set to be $e\big(\, \pi_{m+1} (\psi) \, \big)$ and for $m' > n_k$, $\pi_{m'} (\psi)$ is defined to be $0$.
\end{definition}

The following theorem is proven in \cite{HermoJoosten:2017:TSCArith} and provides a succinct derivability condition between monomial normal forms:

\begin{theorem} 
\label{CharacterizatioMNFDeriv}
For any $\psi_0, \ \psi_1 \in {\sf MNF}$, where $\psi_0 := \MNF{0}{k}$ and $\psi_1 := \la m_0^{\be_0} \ra \top \ \land \ \ldots \ \land \ \la m_j^{\be_j} \ra \top$. We have that $\psi_0 \vdash \psi_1$ iff for any $n < \omega$, $\pi_n (\psi_0) \geq \pi_n (\psi_1)$.
\end{theorem}

\section{A variation on Ignatiev's Frame}

The purpose of this section is to define a modal model $\mathcal{J}$ which is universal for our logic. That is, any derivable sequent will hold everywhere in the model whereas any non-derivable sequent will be refuted somewhere in the model.

The model will be based on specials sequences of ordinals. In order to define them, we need the following central definition.

\begin{definition}
We define \emph{ordinal logarithm} as $\le (0):= 0$ and $\le (\alpha + \omega^\beta) := \beta$.
\end{definition}

With this last definition we are now ready to introduce the set of worlds of our frame.

\begin{definition}
By $\text{Ig}^\omega $ we denote the set of \emph{$\ell$-sequences} or \emph{Ignatiev sequences}. That is, the set of sequences $x :=  \la x_0, x_1, x_2, \ldots \ra$ where for $i \, {<} \, \omega$, $x_{i+1}\leq \le(x_i)$. 
\end{definition}

Given a $\ell$-sequence $x$, if all but finitely many of its elements are zero, we will write $\la x_0 , \ldots, x_n, \vec 0\ra$ to denote such $\le$-sequence or even simply $\la x_0 , \ldots, x_n\ra$ whenever $x_{n+1} =0$.

Next, we can define our frame, which is a minor variation of Ignatiev's frame.

\begin{definition}
$\mathcal{J}_\Lambda := \la I, \{ R_n \}_{n < \omega}\ra$ is defined as follows:
$$I := \{ x \, {\in} \, \text{Ig}^\omega : x_i \, {<} \, \Lambda \ \text{for } i \, {<} \, \omega \}$$ and
$$x R_n y :\Leftrightarrow (\forall \, m \leq n \ x_m \, {>} \, y_m \ \wedge \ \forall \, i \, {>} \, n \ x_i \geq y_i).$$
\end{definition}

Since $\Lambda$ is a fixed ordinal along the paper, from now on we suppress the subindex $\Lambda$. \\

The observations collected in the next lemma all have elementary proofs. Basically, the lemma confirms that the $R_n$ are good to model provability logic and respect the increasing strength of the provability predicates $[n]$.

\begin{lemma}\label{theorem:basicPropertiesRnRelations}\ 
\begin{enumerate}
\item
Each $R_n$ for $n\in \omega$ is
\emph{transitive}: $xR_ny \ yR_nz \ \Rightarrow \ xRz$;

\item\label{theorem:basicPropertiesRnRelations:Item:Noetherian}
Each $R_n$ for $n\in \omega$ is
\emph{Noetherian}: each non-empty $X\subseteq I$ has an $R_n$-maximal element $y\in X$, i.e., $\forall\, x {\in} X\ \neg yR_nx$;

\item
The relations $R_n$ are \emph{monotone} in $n$ in the sense that: $xR_ny \Rightarrow xR_my$ whenever $n>m$.

\end{enumerate}
\end{lemma}

Note that Item (\ref{theorem:basicPropertiesRnRelations:Item:Noetherian}) is equivalent to stating that there are no infinite ascending $R_n$ chains. In other words, the converse of $R_n$ is well-founded. \\
\medskip 

We define the auxiliary relations $R_n^\alpha$ for any $n < \omega$ and $\alpha < \Lambda$. The idea is that the $R_n^\alpha$ will model the $\la n^\alpha \ra$ modality.

\begin{definition} \label{OrdRelation}
Given $x, y \in I$ and $R_n$ on $I$, we recursively define $x R_n^\alpha y$ as follows:
\begin{enumerate}
\item 
$x R_n^0 y \ \ :\Leftrightarrow \ \ x=y$;

\item 
$x R_n^{1+\alpha} y \ :\Leftrightarrow \ \forall \, \beta  {<} 1{+}\alpha  \ \exists z \ \big( xR_nz \ \wedge \ z R_n^\beta y\big)$.
\end{enumerate}
\end{definition}

Let us introduce some simple observations about the $R_n^\al$ relations.

\begin{proposition} Given $x, y \in I$, $n < \omega$ and $\al < \Lambda$:
$$x R_n^{\al + 1} y \ \Leftrightarrow \ \exists z \ \big(xR_nz \ \wedge \ z R_n^\al y \big).$$
\label{SucRelation}
\end{proposition}
\begin{proof}
For the left-to-right implication, assume $x R_n^{\al + 1} y \ $. Therefore, we have that $\forall \, \beta  {<} \alpha{+}1  \ \exists z \ \big( xR_nz \ \wedge \ z R_n^\beta y\big)$, so in particular $\exists z \big(x R_n z \ \wedge \ zR_n^\al y \big)$. For right-to-left implication we proceed analogously. Assume $\exists z \big(x R_n z \ \wedge \ zR_n^\al y \big)$. Thus, $\exists z \ \big(x R_n z \ \wedge \ z R_n^{\al} y \big) \ \wedge \ \forall \be < \al \, \exists z' \ \big( z R_n z' \ \wedge \ z' R_n^\be y \big)$. Hence, we have that $\forall \be < \al + 1 \, \exists z \ \big( x R_n z \ \wedge \ z R_n^{\be} y \big)$, that is, $x R_n^{\al + 1} y$. 
\end{proof}

\begin{proposition} \label{limRelation}
Let $x, y \in I$, $n < \omega$ and $\lambda < \Lambda$ such that $\lambda \in \Lim$:
$$ x R_n^\lambda y \ \Leftrightarrow \ \forall \be < \lambda \ x R_n^{1 + \be} y.$$
\end{proposition}
\begin{proof}
For left-to-right implication, notice that if $x R_n^\lambda y$ then by definition, we have that $\forall \be < \lambda \ \exists u \ \big( x R_n u \ \wedge \ u R_n^{\be} y \big)$. Therefore, in particular,  we obtain that $\forall \be < \lambda \ \exists u \ \big( x R_n u \ \wedge \ u R_n^{1 + \be} y \big)$ thus by transitivity, $\forall \be < \lambda \ \ x R_n^{1 + \be} y $. For the other direction, if $\forall \be < \lambda \ x R_n^{1 + \be} y$,  then in particular, $\forall \be < \lambda \ x R_n^{\be + 1} y$ and then, by Proposition \ref{SucRelation}, $\forall \be < \lambda \ \exists u \ \big( x R_n u \ \wedge \ u R_n^{\be} y \big)$, that is, $x R_n^\lambda y$.
\end{proof}

It is easy to see that for example $\la \omega, \vec 0\ra R^n_0 \la m,\vec 0\ra$ for each $n,m \in \omega$, so that also $\la \omega, \vec 0\ra R^\omega_0 \la m,\vec 0\ra$ for each $m \in \omega$. Clearly, we do not have $\la \omega, \vec 0\ra R^{\omega+1}_0 \la m,\vec 0\ra$ for any $m \in \omega$ but we do have $\la \omega+1, \vec 0\ra R^{\omega+1}_0 \la m,\vec 0\ra$ for all $m \in \omega$.

We also note that the dual definition $x \overline R_n^0 y \ :\Leftrightarrow \ x=y$; and $x \overline R_n^{1+\alpha} y \ :\Leftrightarrow \ \forall \, \beta  {<} 1{+}\alpha  \ \exists z \ \big( x\overline R^\beta_nz \ \wedge \ z \overline R_n y\big)$ does not make much sense on our frames. For example we could have $\la \omega,\vec 0\ra \overline R^{\alpha}_0 \la 0,\vec 0\ra$ for any ordinal $\alpha >0$.

With the the auxiliary relations $R_n^\alpha$, we give the following definition for a formula $\varphi$ being true in a point $x$ of $\mathcal{J}$. 

\begin{definition} \label{FormulaTrueInPoint}
Let $x \in I$ and $\varphi \in \FLE$. By $x \Vdash \varphi$ we denote the validity of $\varphi$ in $x$ that is recursively defined as follows:
\begin{itemize}
\item $x \Vdash \top$ for all $x \in I$; 
\item $x \Vdash \varphi \wedge \psi$ iff $x \Vdash \varphi$ and $x \Vdash \psi$;
\item $x \Vdash \la n^\alpha \ra \varphi$ iff there is $y \in I$, $x R_n^\alpha y$ and $y \Vdash \varphi$.
\end{itemize}
\end{definition}
 
Here are some easy observations on the $R^\alpha_n$ relations which among others tell us that all the $R^\alpha_n$ serve the purpose of a provability predicate for any $n\in \omega$ and $\alpha < \Lambda$.

\clearpage
\thispagestyle{empty}
\begin{figure}
\caption{A fragment of $\mathcal{J}$. The dashed arrows represent $R_0$ relations, while the continuous arrows represent $R_1$ relations.}
\hspace{-3cm}\includegraphics{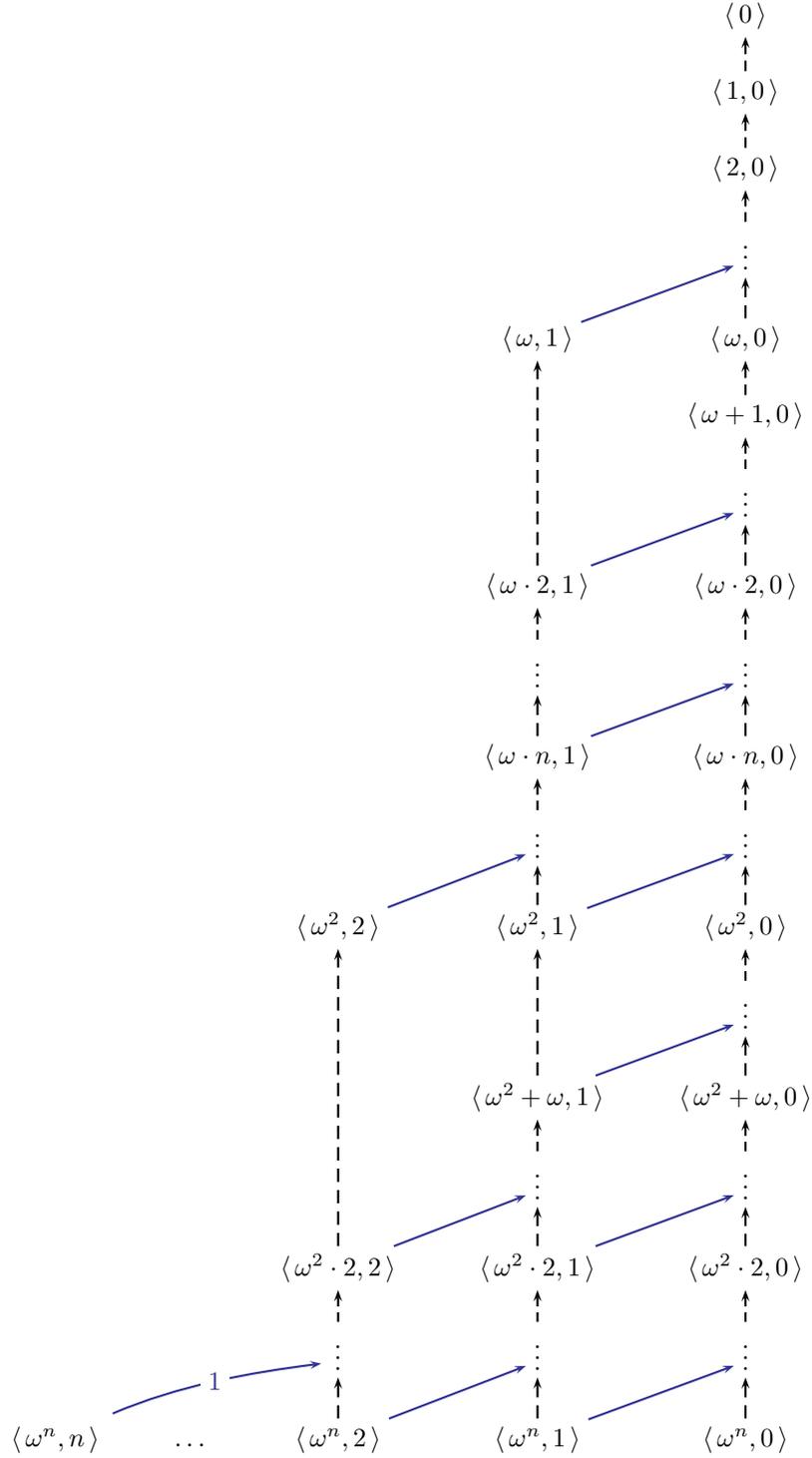}
\end{figure}

\clearpage
\setcounter{page}{7}.

\begin{lemma}\label{theorem:basicPropertiesRnAlphaRelations}\ 
\begin{enumerate}
\item
Each $R_n^{1+\alpha}$ for $n\in \omega$ and ${\alpha}$ an ordinal is
\emph{transitive}: $xR_n^{1+\alpha} y \ \wedge \  yR^{1+\alpha}_nz \ \Rightarrow \ xR_n^{1+\alpha} z$;

\item\label{theorem:basicPropertiesRnRelations:Item:Noetherian}
Each $R_n^{1+\alpha}$ for $n\in \omega$ and ${\alpha}$ an ordinal is
\emph{Noetherian}: each non-empty $X\subseteq I$ has an $R^{1+\alpha}_n$-maximal element $y\in X$, i.e., $\forall\, x {\in} X\ \neg yR_n^{1+\alpha} x$;

\item
The relations $R_n^{1+\alpha}$ are \emph{monotone} in $n$ in the sense that: $xR_n^{1+\alpha} y \Rightarrow xR_m^{1+\alpha} y$ whenever $n>m$;

\item \label{theorem:basicPropertiesRnRelations:Item:Mon}

The relations $R_n^{1+\alpha}$ are \emph{monotone} in ${1+\alpha}$ in the sense that: $xR_n^{1+\alpha} y \Rightarrow xR_n^{1+\beta} y$ whenever $1+ \beta < {1+\alpha}$.
\end{enumerate}
\end{lemma}

\begin{proof}
The first three items follow directly from Lemma \ref{theorem:basicPropertiesRnRelations} by an easy transfinite induction. The last item is also easy.
\end{proof}

\section{A characterization for transfinite accessibility}

The intuitive idea between the $xR^\alpha_n y$ assertion, is that this tells us that there exists a chain of `length' $\alpha$ of $R_n$ steps leading from the point $x$ up to the point $y$. The following useful lemma tries to capture this intuition. 

\begin{lemma}\label{theorem:RalphaCharacterizationViaChains}
For $x,y \in I$ and $n<\omega$ we have that the following are equivalent
\begin{enumerate} 
\item
$xR^{1+\alpha}_n y$

\item
For each $\beta< 1+\alpha$ there exists a collection $\{ x^\gamma\}_{\gamma < \beta}$ so that 
\begin{enumerate}
\item
$xR_n x^{\gamma}$ for any $\gamma < \beta$,

\item
$x^0 = y$ and,

\item
for any $\gamma' < \gamma < \beta$ we have $x^\gamma R_n x^{\gamma'}$.
\end{enumerate}
\end{enumerate}
\end{lemma}

\begin{proof}
By induction on $\alpha$.
\end{proof}

We shall now provide a characterization of the $R^{1+\alpha}_n$ relations. To this end, let us for convenience define 
\[
x R_{-1}^{\zeta} y \ \ :\Leftrightarrow \ \ \forall \, n{>}0 \ \ x_n \geq y_n.
\]
With this notation the following theorem makes sense.

\begin{theorem}\label{theorem:CharacterizationOfRnAlphaSteps}
For $x,y \in I$ and $n<\omega$ we have that the following are equivalent
\begin{enumerate} 
\item \label{Char1}
$ xR^{1+\alpha}_n y$;

\item \label{Char2}
$ x_n \geq y_n + \big(1+ e(y_{n+1})\big)\cdot(1+\alpha)$ and $xR_{n-1}^{e(1+\alpha)}y$;

\item \label{Char3}
\[
\begin{array}{ll}
x_n \geq y_n + \big(1+ e(y_{n+1})\big)\cdot(1+\alpha) &\mbox{ and,}\\
\mbox{$x_m>y_m$ for $m<n$} & \mbox{ and,}\\
\mbox{$x_m\geq y_m$ for $m>n$.} & \\
\end{array}
\]
\end{enumerate}
\end{theorem}

We dedicate the remainder of this section to proving this theorem and move there through a series of lemmas. The first lemma in this series is pretty obvious. It tells us that if we can move from $x$ to $y$ in $\alpha$ many steps, then the distance between $x_n$ and $y_n$ must allow $\alpha$ many steps; That is, they lie at least $\alpha$ apart.

\begin{lemma}\label{theorem:RnAlphaStepsImplyTheNcoordinateToLieAlphaApart}
For $x,y \in I$ and $n<\omega$ and any ordinal $\alpha < \Lambda$, if $xR_n^\alpha y$ then $ x_n \geq y_n + \alpha$.
\end{lemma}

\begin{proof}
By an easy induction on $\alpha$.
\end{proof}

However, how many $R_n$ steps one can make is not entirely determined by the $n$ coordinates of the points. For example, there is just a single $R_0$ step from the point $\la \omega \cdot 2, 1 \ra$ to the point $\la \omega, 1 \ra$ whereas these points lie $\omega$ apart on the `$0$ coordinate'. The following lemma tells us how for $R_n$ steps, the $n$-th coordinates are affected by the values of the $n+1$-th coordinate.

\begin{lemma}\label{theorem:manyNsteps}
For $x,y \in I$ and $n<\omega$ with $xR^{1+\alpha}_ny$, we have
\[
x_n \geq y_n + e(y_{n+1})\cdot(1+\alpha).
\]
\end{lemma}

In order to give a smooth presentation of this proof, we first give two simple technical lemmas with useful observations on the ordinals and ordinal functions involved.

\begin{lemma}\label{theorem:TechnicalLemmaOnOrdinals}
For $\alpha, \beta$ and $\gamma$ ordinals we have
\begin{enumerate}
\item
$\le (\beta) \geq 1+ \alpha \ \ \Longleftrightarrow \ \ \beta \in e(1+ \alpha) \cdot (1 + \On)$, \label{techOrd1}

\item
If $(1+\alpha)< \beta$ and $\gamma \in e(\beta)\cdot (1 + \On)$, then $\gamma \in e(1+\alpha)\cdot (1 + \On)$, \label{techOrd2}

\item
$e(\beta + (1+\alpha)) = e(\beta)\cdot e(1+ \alpha)$, \label{techOrd3}

\item \label{techOrd4}
For $\alpha$ a limit ordinal, we have that 
\[
xR^\alpha_n y \ \ \Longleftrightarrow \ \ \forall \, 1{+}\beta {<} \alpha\, \exists z \ (xR_nz \wedge zR_n^{1+\beta}y).
\]
\end{enumerate}
\end{lemma}

\begin{proof}
The first two items then can easily be seen by using a Cantor Normal Form expression with base $\omega$. For Item (\ref{techOrd1}), we use the fact that $\beta \in \Lim$ together with that if $\le (\beta) \geq 1+ \alpha$ , then $\beta \geq e(\le(\beta)) \geq e(1+ \alpha)$. For Items (\ref{techOrd2}) and (\ref{techOrd3}) we use that $e(1+\omega) = \omega^{1+\omega} = \omega^1 \cdot \omega^\omega$. The last item follows from Definition \ref{OrdRelation} together with the fact that $1 + \alpha \in \Lim = \alpha$.
\end{proof}

\begin{lemma}
For $x, y \in I$ and $n < \omega$, $xR_n y \Longleftrightarrow x_n \geq y_n + e(x_{n+1})$. 
\end{lemma}
\begin{proof}
We make a case distinction on $x_n$. If $x_n \in \text{Succ}$ then is trivial since $e(x_{n+1}) = 0$. If $x_n \in \Lim$, and furthermore, $x_n$ is an additively indecomposable limit ordinal, it follows from the fact that $x_n > y_n$ and $x_n \geq e(x_{n+1})$. Otherwise, we can rewrite $x_n$ as $\alpha + e(\beta)$ for some $\beta \geq x_{n+1}$, and $y_n$ as $\delta + \omega^\gamma$.  If $y_n \leq \alpha$ then clearly $x_n \geq y_n + e(x_{n+1})$. If $\alpha = \delta$ and $\gamma <  \beta$, then notice that $\omega^\gamma + e(\beta) = e(\beta)$ Thus, we have that $\alpha + e(\beta) = \delta + \omega^\gamma + e(\beta) \geq y_n + e(x_{n+1})$.

\end{proof}

With these technical lemmas at hand we can now prove Lemma \ref{theorem:manyNsteps}.

\begin{proof}
By induction on $\alpha$. For $\alpha := 0$, we check that $x_n \geq y_n + e(y_{n+1})$. Note that since $x R_n y$ then $x_n \geq y_n + e(x_{n+1})$ and $x_{n+1} \geq y_{n+1}$, then $x_n \geq y_n + e(y_{n+1})$. For $\alpha := \beta + 1$, if $x R_n^{1+\beta} y$ then there is $z \in I$ such that $x R_n z$ and $z R_n^{1 + \beta} y$. Thus, we have the following:
\begin{enumerate}
\item $x_n \geq z_n + e(z_{n+1})$;
\item $z_n \geq y_n + e(y_{n+1}) \cdot (1 + \beta)$.
\end{enumerate}
Therefore, $x_n \geq y_n + e(y_{n+1}) \cdot (1 + \beta) + e(z_{n+1})$. Since $e(z_{n+1}) \geq e(y_{n+1})$ then $x_n \geq y_n + e(y_{n+1}) \cdot (1 + \beta) + e(y_{n+1})$ i.e.  $x_n \geq y_n + e(y_{n+1}) \cdot (1 + \beta + 1)$. For $\alpha \in \Lim$, notice that by IH, we have that $x_n \geq y_n + e(y_{n+1}) \cdot (1 + \delta)$ for $\delta < \alpha$. Thus, $x_n \geq y_n + e(y_{n+1}) \cdot (1 + \alpha)$.
\end{proof}

Combining Lemma \ref{theorem:manyNsteps} and Lemma \ref{theorem:RnAlphaStepsImplyTheNcoordinateToLieAlphaApart} we get the following.

\begin{corollary}\label{theorem:AlphaNstepImpliesAbigNcoordinat}
For $x,y \in I$ and $n<\omega$ we have that 
\[
xR^{1+\alpha}_n y \Rightarrow x_n \geq y_n + \big(1+ e(y_{n+1})\big)\cdot(1+\alpha).
\]
\end{corollary}

This corollary takes care of part of the implication from Item (\ref{Char1}) to Item (\ref{Char2}) in Theorem \ref{theorem:CharacterizationOfRnAlphaSteps}. We will no focus on the implication from Item (\ref{Char3}) to Item (\ref{Char1}) but before we do so, we first formulate a simple yet useful lemma.

\begin{lemma}\label{theorem:nextCoordinatesGetSubstantiallyBigger}
For $x,y \in I$, if $xR_{m+1}y$, then $x_m \geq y_m + e(x_{m+1})$.
\end{lemma}

\begin{proof}
Since $xR_{m+1}y$, in particular $x_{m+1}>0$ whence $x_m \in e(x_{m+1}) \cdot (1 + {\sf On})$ and the result follows by writing both $x_m$ and $y_m$ in Cantor Normal Form.
\end{proof} 

With this technical lemma we can obtain the next step in the direction from Item (\ref{Char3}) to Item (\ref{Char1}) in Theorem \ref{theorem:CharacterizationOfRnAlphaSteps}.

\begin{lemma}\label{theorem:TheNcoordinatesFarApartImplyManyNsteps}
For $x,y \in I$ and $n<\omega$ we have that if
\[
\begin{array}{ll}
x_n \geq y_n + \big(1+ e(y_{n+1})\big)\cdot(1+\alpha) &\mbox{ and,}\\
\mbox{$x_m>y_m$ for $m<n$} & \mbox{ and,}\\
\mbox{$x_m\geq y_m$ for $m>n$.} & \\
\end{array}
\]
then
\[
xR^{1+\alpha}_n y.
\]
\end{lemma}

\begin{proof}
We use Lemma \ref{theorem:RalphaCharacterizationViaChains} whence are done if we can find for each $\beta< 1+\alpha$ there exists a collection $\{ x^\gamma\}_{\gamma < \beta}$ so that 
\begin{enumerate}
\item
$xR_n x^{\gamma}$ for any $\gamma < \beta$,

\item
$x_0 = y$ and,

\item
for any $\gamma' < \gamma < \beta$ we have $x^\gamma R_n x^{\gamma'}$.
\end{enumerate}
We define $x^\gamma$ uniformly as follows. We define $x^0 := y$ and
\[
x^{1+\gamma}_m \ \ : = \ \   
\begin{cases}  y_m & \text{in case $m>n$,} \\ 
y_m + \big(1+ e(y_{n+1})\big)\cdot (1+ \gamma) & \text{in case $m=n$,} \\  
y_m + e(y_{m+1}) &\text{in case $m<n$.} \\ 
\end{cases} 
\]
We make a collection of simple observations:
\begin{enumerate}[i]
\item
Each $x^\gamma$ is an element of $I$ for any $\gamma < \alpha$ since $x^\gamma_{m+1} \leq \le (x^\gamma_m)$ for any $m$;

\item
We now see that $xR_nx^\gamma$ for each $\gamma < \alpha$. For $m > n$ we obviously have that $x_m \geq x^\gamma_m$ and also $x_n > x^\gamma_n$ is clear. By induction we see that $x_m > x_m^\gamma$ using Lemma \ref{theorem:nextCoordinatesGetSubstantiallyBigger} and the fact that $e$ is a strictly monotonously growing ordinal function;

\item
$x_0 =y$ by definition;

\item
By strict monotonicity of $e$, we see that for any $\gamma' < \gamma < \alpha$ we have $x^\gamma R_n x^{\gamma'}$.

\end{enumerate} 
\end{proof}

We are now ready to prove Theorem \ref{theorem:CharacterizationOfRnAlphaSteps}.

\begin{proof}
From Item (\ref{Char2}) to Item (\ref{Char3}) is easy and from Item (\ref{Char3}) to Item (\ref{Char1}) is Lemma \ref{theorem:TheNcoordinatesFarApartImplyManyNsteps} so we focus on the remaining implication.

As mentioned before, half of the implication from Item (\ref{Char1}) to Item (\ref{Char2}) follows from Corollary \ref{theorem:AlphaNstepImpliesAbigNcoordinat} so that it remains to show that $xR^{1+\alpha}_ny \ \Rightarrow \ xR^{e(1+\alpha)}_{n-1}y$. For $n=0$ this is trivial and in case $n\neq 0$ we reason as follows. 

Since $xR^{1+\alpha}_n y$ we get in particular that $x_n \geq y_n + 1 + \alpha$. Thus, by Lemma \ref{theorem:nextCoordinatesGetSubstantiallyBigger} we see 
\[
x_{n-1} \geq y_{n-1} + e(x_n) \geq y_{n-1} + e(y_n + 1 + \alpha).
\]
Now using the fact (Lemma \ref{theorem:TechnicalLemmaOnOrdinals}) that $e(y_n + 1 + \alpha) = e(y_n)\cdot e (1+ \alpha)$ we see, making a case distinction whether $y_n =0$ or not and using that $e(1 + \alpha)$ is a limit ordinal, that 
\[
x_{n-1} \geq y_{n-1} + (1+e(y_n)\cdot ( 1+  e(1 + \alpha)).
\]
The result now follows from an application of Lemma \ref{theorem:TheNcoordinatesFarApartImplyManyNsteps}.
\end{proof}

\section{Definable sets}

In this section we shall define a translation between formulas in {\sf MNF} and Ignatiev sequences with finite support as well as a way of characterizing subsets of $I$. Moreover, we shall see how some of these subsets of $I$ can be related to the extensions of formulas.

\begin{definition}\label{TranslationMNFI}
Let $\psi:= \MNF{0}{k} \in {\sf MNF}$. By $x_\psi$ we denote the sequence $\la \pi_i (\psi) \ra_{i < \omega}$.
\end{definition}

In virtue of Definition \ref{projection}, we can observe that for $\psi \in {\sf MNF}$, we have that $x_\psi \in \text{Ig}^\omega $. Furthermore, we shall see that $x_\varphi$ is the ``first'' point in $I$ where $\varphi$ holds. First we can make some simple observations.

\begin{lemma} \
\begin{enumerate}
\item For any $x \in I$, $x \Vdash \la n^\al \ra \top$ iff $x_n \geq \al$; \label{trueInTrans1}
\item For any $\psi \in {\sf MNF}$, $x_\psi  \Vdash \psi$. \label{trueInTrans2}
\end{enumerate}
\label{trueInTranslation}
\end{lemma}
\begin{proof}
The second item follows from the first one and Definition \ref{TranslationMNFI}. For the right-to-left implication of the first item, assume $x_n \geq \al \, {>} \, 0$. Therefore, for $i \, {<} \, n$, we have that $x_i \, {>} \, 0$ and for $i' \, {>} \, n$, $x_i \geq 0$. Thus, by Theorem \ref{theorem:CharacterizationOfRnAlphaSteps}, $x R_n^\al \langle 0 \rangle$ and so $x \Vdash \la n^\al \ra \top$. For the other direction, assume $x \Vdash \la n^\al \ra \top$ for $\al \, {>} \, 0$. Hence, there is $y \in I$ such that $x R_n^\al y$ and $y \Vdash \top$. By Theorem \ref{theorem:CharacterizationOfRnAlphaSteps}, $x_n \geq y_n + \big( 1 + e(y_{n+1}) \big) \cdot \al$ and so, $x_n \geq \al$. The case $\al = 0$ is straightforward.  
\end{proof}

The following two definitions introduce the extension of Ignatiev sequences and the extension of formulas, respectively.

\begin{definition}
Given $x \in I$, by $\llbracket x \rrbracket$ we denote the set of $\ell$-sequences which are coordinate-wise at least as big as $x$. That is, we define $\llbracket x \rrbracket:= \{ y \,{\in}\, I : y_i \geq x_i \text{ for every } i\, {<} \, \omega \}$. \label{extensionWorld}
\end{definition}

\begin{definition}
Let $\varphi \in \FLE$. By $\llbracket\, \varphi \, \rrbracket$ we denote the set of worlds where $\varphi$ holds i.e.  $\llbracket\, \varphi \, \rrbracket = \{x \in I : x \Vdash \varphi \}$. \label{extensionFormula}
\end{definition}

The following lemma relates definitions \ref{extensionWorld} and \ref{extensionFormula}.   

\begin{lemma}
For any $\varphi \in \FLE$, there is $x := \langle x_0,\, \ldots\, ,x_k, 0 \rangle \in I$ such that \linebreak $\llbracket\, \varphi \, \rrbracket = \llbracket x \rrbracket$. \label{extensionformulaWorld}
\end{lemma}
\begin{proof}
The proof goes by induction on $\varphi$. The base case is trivial. For \textbf{the conjunctive case}, let $\varphi = \psi \land  \chi$. By the I.H. we have that there are $y, \, z \in I$ such that $\llbracket\, \psi \, \rrbracket = \llbracket y \rrbracket$ and $\llbracket\, \chi \, \rrbracket = \llbracket z \rrbracket$. Moreover, by the I.H. we also have that $y := \langle y_0,\, \ldots\, ,y_j, 0 \rangle $ and $z := \langle z_0,\, \ldots\, ,z_i, 0 \rangle$. Let $n$ be the index of the rightmost non-zero component. Hence we can define $x$ as follows:
\begin{itemize}
\item $x_i = \max(y_i, z_i)$ for $i \geq n$;
\item $x_{i} = \min\{ \delta : \delta \geq \max(y_{i}, z_{i}) \ \& \ l(\delta) \geq x_{i+1}\}$ for $i \, {<} \, n$.

\end{itemize} 
We can easily check that $x \in I$. Next, we check that for any $x' \in I$, we have that $x' \Vdash \psi \land \chi$ iff $x' \in \llbracket x \rrbracket$. For right-to-left implication, consider $x' \in \llbracket x \rrbracket$. Thus, for $k \, {<} \, \omega$, we have that both $x'_k \geq x_k \geq y_k$ and $x'_k \geq x_k \geq z_k$. Thus, $x' \in \llbracket y \rrbracket \,  \cap \, \llbracket z \rrbracket $ and so by the I.H. $x' \Vdash \psi \land \chi$. For the other direction, consider $x' \in I$ such that $x' \Vdash \psi \land \chi$. Clearly, for $i > n$, we have that $x'_i \geq x_i$. We check by induction on $k$ that $x'_{n-k} \geq x_{n-k}$. For the base case, since $x'\Vdash \psi \land \chi$, then by the I.H. $x' \in \llbracket y \rrbracket \,  \cap \, \llbracket z \rrbracket $ and so  $x'_n \geq y_n$ and $x'_n \geq z_n$. Thus,  $x'_n \geq \max(y_n, z_n) = x_n$. For the inductive step, by definition of Ignatiev sequences together with the I.H., we have that $l(x'_{n-(k + 1)}) \geq x'_{n- k} \geq x_{n- k}$ and since $x'\Vdash \psi \land \chi$, then $x'_{n-(k+1)} \geq \max(y_{n-(k+1)}, z_{n-(k+1)})$. Therefore, being $x_{n- (k+1)}$ the minimal ordinal satisfying both conditions, we can conclude that $x'_{n- (k+1)} \geq x_{n- (k+1)}$. Hence, $\llbracket\, \psi \land \chi \, \rrbracket = \llbracket x \rrbracket$.

For \textbf{the modality case}, let $\varphi := \la n^\al \ra \psi$ with $\al\, {>} \, 0$. Thus, by the I.H. there is $y \in I$ such that $\llbracket \, \psi \, \rrbracket = \llbracket y \rrbracket$ and $y := \langle y_0,\, \ldots\, ,y_j, 0 \rangle$. We can define $x$ as follows:
\begin{itemize}
\item $x_i = y_i$ for $i\, {>} \, n$;
\item $x_n = y_n + \big(1 + e(y_{n+1}) \big) \cdot \al$;
\item $x_{i} = \min \{ \delta : \delta \geq y_{i} \ \& \ l(\delta) \geq x_{i+1} \}$ for $i \, {<} \, n$.  
\end{itemize}
As in the previous case, we can easily check that $x \in I$. We claim that $\llbracket x \rrbracket = \llbracket \la n^\al \ra \psi \rrbracket$. Let $x' \in \llbracket x \rrbracket$. By Theorem \ref{theorem:CharacterizationOfRnAlphaSteps} we can see that $x R_n^\al y$. Hence, since $x'_i \geq x_i$ for $i \, {<} \, \omega$, $x' R_n^\al y$ and so $x' \Vdash \la n^\al \ra \psi$. For the other inclusion, consider $x' \in I$ such that $x' \Vdash \la n^\al \ra \psi$. By the I.H. and Theorem \ref{theorem:CharacterizationOfRnAlphaSteps}, we can easily check that for $i \, {>} \, n$, we have that $x'_i \geq x_i$. For $i \leq n$, we proceed by an easy induction on $k$ to see that $z_{n-k} \geq x_{n-k}$. The base case follows directly from Theorem \ref{theorem:CharacterizationOfRnAlphaSteps}. For the inductive step, by definition of Ignatiev sequences together with the I.H., we have that $l(x'_{n-(k + 1)}) \geq x'_{n- k} \geq x_{n- k}$. Since $x' \Vdash \la n^\al \ra \psi$, then there is $z \in I$ such that $x R_n^\al z$ and $z \Vdash \psi$. Thus, by the I.H., $z \in \llbracket y \rrbracket$, and so we have that $x'_{n-(k+1)} \, {>} \, z_{n- (k+1)} \geq y_{n-(k+1)}$. Therefore, we get that $l(x'_{n-(k + 1)}) \geq x_{n- k}$ and $x'_{n-(k+1)} \, {>} \,  y_{n-(k+1)}$. Thus, since $x_{n- (k+1)}$ is the least ordinal satisfying both conditions, we have that $x'_{n- (k+1)} \geq x_{n- (k+1)}$. 
\end{proof}

\section{Soundness}

To prove the soundness of $\TPr$, let us begin by semantically define the entailment between our modal formulas. 

\begin{definition}
For any formulas $\varphi, \, \psi \in \FLE$, we write $\varphi \models \psi$  iff for all $x \in I$, if $x \Vdash \varphi$ then $x \Vdash \psi$. Analogously, we write $\varphi \equiv_{\mathcal{I}} \psi$ iff for any $x \in I$, we have that $x \Vdash \varphi$ iff $x \Vdash \psi$.  
\end{definition}

With our notion of semantical entailment we can formulate our soundness theorem.

\begin{theorem}[Soundness] \label{Soundness}
For any formulas $\varphi, \, \psi \in \FLE$, if $\varphi \vdash \psi$ then \linebreak $\varphi \models \psi$.
\end{theorem}

\begin{proof}
By induction on the length of a $\TPr$ proof of $\varphi \vdash \psi$. It is easy to see that the first three rules preserve validity. With respect to the axioms, the first two axioms are easily seen to be valid. The the correctness of reduction axiom is given by Theorem \ref{theorem:CharacterizationOfRnAlphaSteps}. The remaining axioms and rules are separately proven to be sound in the remainder of this section.
\end{proof}

We start by proving the soundness of co-additivity axiom i.e.
\[
\la n^\al \ra \la n^\be \ra \varphi \equiv_\mathcal{I} \la n^{\be + \al} \ra \varphi.
\]

\begin{proposition}
For any $x, y, z \in I$, $n < \omega$ and $\al, \be < \Lambda$ we have that $x R_n^\al y$ and $y R_n^\be z$ iff $x R_n^{\be + \al} z$. \label{soundCoadd}
\end{proposition}
\begin{proof}
We proceed by transfinite induction on $\al$ with the base case being trivial. For $\al \in \suc$, let $\al := \delta + 1$ for some $\delta$. Therefore:
\begin{tabbing}
$x R_n^\al y$ and $y R_n^\be z \ $ \= $\Leftrightarrow \ x R_n^{\delta + 1} y$ and $y R_n^\be z$; \\[0.30cm]
\ \> $\Leftrightarrow \ \exists u \ \big( x R_n u \ \wedge \ u R_n^\delta y \ \wedge \ y R_n^\be z \big)$; \\[0.30cm]
\ \> $\Leftrightarrow \ \exists u \ \big( x R_n u \ \wedge \  u R_n^{\be + \delta} z \big)$, by the I.H. ; \\[0.30cm]
\ \> $\Leftrightarrow \ x R_n^{\be + \delta + 1} z $; \\ [0.30cm]
\ \> $\Leftrightarrow \ x R_n^{\be + \al} z $.
\end{tabbing}
For $\al \in \Lim$, we have that $x R_n^\al y$ and $y R_n^\be z \ \Leftrightarrow \ \forall \delta < \al \ \big( x R_n^{1 + \delta} y \ \wedge \ y R_n^\be z \big)$ by Proposition \ref{limRelation}. By the I.H. we obtain $\forall \delta < \al \ x R_n^{\be + 1 + \delta} z$ and so $x R_n^{\be + \al} z$. 
\end{proof}

With this last result, we get the co-additivity of the $R_n^\al$ relations. This together with Definition \ref{FormulaTrueInPoint} gives us the following corollary.

\begin{corollary}
The co-additivity axiom is sound.
\end{corollary}
\begin{proof}
By Definition \ref{FormulaTrueInPoint}, $x \Vdash \la n^\al \ra \la n^\be \ra \varphi$ iff there are $y, \, z \in I$ such that $x R_n^\al y$, $y R_n^\be z$ and $z \Vdash \varphi$. Thus, by Proposition \ref{soundCoadd}, $x \Vdash \la n^\al \ra \la n^\be \ra \varphi$ iff $x R_n^{\be + \al} z$ and $z \Vdash \varphi$ i.e. $x \Vdash \la n^{\be + \al} \ra \varphi$. \end{proof}

\begin{proposition} The monotonicity axiom is sound, that is: 
\[
\la n^\al \ra \varphi \models \la n^\be \ra \varphi
\]
for $\be \, {<} \, \al$.
\end{proposition}
\begin{proof}
With the help Lemma \ref{theorem:basicPropertiesRnAlphaRelations}, Item (\ref{theorem:basicPropertiesRnRelations:Item:Mon}), we have that if $x \Vdash \la n^\al \ra \varphi$ then $x \Vdash \la n^\be \ra \varphi$ for $\be, \ 0 < \be < \al$. We check that if $x \Vdash \la n^1 \ra \varphi$ then  $x \Vdash \varphi$ by induction on $\varphi$. 

The Base and the conjunctive cases are straightforward, so we consider $\varphi := \la m^\delta \ra \psi$ and assume $x \Vdash \la n^1 \ra \la m^\delta \ra \psi$. We make the following case distinction:
\begin{itemize}
\item If $n = m$, then by soundness of co-additivity axiom together with Lemma \ref{theorem:basicPropertiesRnAlphaRelations}, Item (\ref{theorem:basicPropertiesRnRelations:Item:Mon}) we have that $x \Vdash \la m^\delta \ra \psi$;
\item If $n > m$, then monotonicity property of $R_n^{1 + \al}$ together with soundness of co-additivity axiom and Lemma \ref{theorem:basicPropertiesRnAlphaRelations}, Item (\ref{theorem:basicPropertiesRnRelations:Item:Mon}) we have that $x \Vdash \la m^\delta \ra \psi$;
\item If $n < m$, then there are $y, \, z \in I$ such that $x\,  R_n \, y \, R_m^\delta\, z$ and $z \Vdash \psi$. Thus, we can easily check that $x \, R_m^\delta \, z$, and so $x \Vdash \la m^\delta \ra \psi$.

\end{itemize}
\end{proof}

The following proposition establishes the correction of the Schmerl axiom by using the translation between formulas in monomial normal form and Ignatiev sequences.

\begin{proposition} The Schmerl axiom is sound i.e.
\[ 
\la n^\al \ra \big( \, \la n_0^{\al_0} \ra \top \ \wedge \ \psi \, \big) \equiv_{\mathcal{I}} \la n^{e^{n_0 - n} (\al_0) \cdot (1 + \al)}  \ra \top \land \la n_0^{\al_0} \ra \top \ \land \ \psi
\]
for $n < n_0$ and $\la n_0^{\al_0} \ra \top \ \wedge \ \psi \in {\sf MNF}$.
\end{proposition}
\begin{proof}
For the left-to-right direction, assume $x \Vdash \la n^\al \ra \big( \la n_0^{\al_0} \ra \top \ \land \ \psi \big)$. Thus, by soundness of monotonicity axiom, we have that $x \Vdash \la n_0^{\al_0} \ra \top \ \land \ \psi$. Therefore, we only need to check that $x \Vdash \la n^{e^{n_0 - n}(\al_0) \cdot (1+ \al)} \ra \top$. Notice that $x \Vdash \la n^\al \ra \la n_0^{\al_0} \ra \top$ and so there are $y,\, z \in I$ such that $x R_n^\al y R_{n_0}^{\al_0} z$. By Theorem \ref{theorem:CharacterizationOfRnAlphaSteps} we have that 
\begin{equation}
x_n \geq y_n + (1 + e(y_{n+1})) \cdot \al. \label{ineqAx6}
\end{equation} 
Also notice that since $y R_{n_0}^{\al_0} z$ then $y R_n^{e^{n_0 - n} (\al_0)} z$ and $y R_{n+1}^{e^{n_0 - {n+1}} (\al_0)} z$. Hence by Theorem \ref{theorem:CharacterizationOfRnAlphaSteps} $y_n \geq e^{n_0 - n} (\al_0)$ and $y_{n+1} \geq e^{n_0 - {n+1}} (\al_0)$. Combining this with \ref{ineqAx6} we get that $x_n \geq e^{n_0 - n} (\al_0) + \big(1 + e(e^{n_0 - {n+1}}(\al_0))\big) \cdot \al = e^{n_0 - n} (\al_0) \cdot (1 +\al)$. Thus, in particular, we have that $x R_n^{e^{n_0 - n} (\al_0) \cdot (1 +\al)} \langle 0 \rangle$ and so, $x \Vdash \la n^{e^{n_0 - n}(\al_0) \cdot (1+\al)} \ra \top$. \\

For the other direction, assume $x \Vdash \la n^{e^{n_0 - n}(\al_0) \cdot (1+\al)} \ra \top  \land  \la n_0^{\al_0} \ra \top  \land  \psi$. Hence, $x \Vdash \la n^{e^{n_0 - n}(\al_0) \cdot (1+\al)} \ra \top$ and so, by Lemma \ref{trueInTranslation}, Item (\ref{trueInTrans1}), $x_n \geq e^{n_0 - n}(\al_0) \cdot (1+\al) = e^{n_0 - n}(\al_0) +(1 + e^{n_0 - n}(\al_0)) \cdot \al$. Since $\la n_0^{\al_0} \ra \top  \land  \psi\in {\sf MNF}$ consider $y_{\la n_0^{\al_0} \ra \top  \land  \psi}$. Notice that $\pi_n (\la n_0^{\al_0} \ra \top  \land  \psi) = e^{n_0 - n} (\al_0)$, thus by Defintion \ref{TranslationMNFI} and Theorem \ref{theorem:CharacterizationOfRnAlphaSteps} we can easly check that $x R_n^\al y_{\la n_0^{\al_0} \ra \top  \land  \psi}$ and by Lemma \ref{trueInTranslation}, Item (\ref{trueInTrans2}), $y_{\la n_0^{\al_0} \ra \top  \land  \psi} \Vdash \la n_0^{\al_0} \ra \top  \land  \psi$. Therefore, $x \Vdash \la n^\al \ra \big( \, \la n_0^{\al_0} \ra \top \ \wedge \ \psi \, \big)$.  
\end{proof}

Lastly, we check the soundness of Rule (\ref{r:4}) by applying the relation between definable sets and the extension of Ignatiev sequences proved in Lemma \ref{extensionformulaWorld}. This next result concludes the soundness proof of $\TPr$.

\begin{proposition} If $\varphi \models \psi$ then, for $m < n$: 
\[
\la n^\al \ra \varphi \, \land \, \la m^{\be + 1} \ra \psi \, \models \, \la n^\al \ra \big( \, \varphi \, \land \, \la m^{\be+ 1} \ra \psi \, \big).
\]
\end{proposition}
\begin{proof}
Assume $\varphi \models \psi$ and let $x \in I$ such that $x \Vdash \la n^\al \ra \varphi \, \wedge \, \la m^{\be + 1} \ra \psi$. Since $\varphi \models \psi$, by Lemma \ref{extensionformulaWorld}, there are $y, \, z \in I$ such that $\llbracket y \rrbracket = \llbracket \varphi \rrbracket \subseteq \llbracket \psi \rrbracket = \llbracket z \rrbracket$. Let $y', \, z' \in I$ such that $\llbracket y' \rrbracket = \llbracket \la n^\al \ra \varphi \rrbracket$ and $\llbracket z' \rrbracket = \llbracket \la m^{\be + 1} \ra \psi \rrbracket$, and $w \in I$ such that $\llbracket w \rrbracket = \llbracket  \varphi \, \land \, \la m^{\be+ 1} \ra \psi \rrbracket$. Since $y \in \llbracket z \rrbracket$, we know that $w_i =y_i$ for $i \, {>} \, m$. For the remaining components, we have that:
\begin{itemize}
\item $w_m = \max\big(y_m, z'_m \big) $;
\item $w_i = \min\{ \delta : \delta \geq \max(y_i, z'_i) \ \& \ l(\delta) \geq w_{i+1}\}$ for $i \, {<} \, m$.
\end{itemize} 
On the other hand, since $x \Vdash \la n^\al \ra \varphi \, \wedge \, \la m^{\be + 1} \ra \psi$, we have the following:
\begin{itemize}
\item $x_i \geq y_i$ for $i \, {>} \, n$;
\item $x_n \geq y'_n$;
\item $x_i \geq \min \{ \delta : \delta \geq y'_i \ \& \ l(\delta) \geq x_{i+1} \}$ for $i, \ m \, {<} \, i \, {<} \, n$;
\item $x_i \geq \min\{ \delta : \delta \geq \max(y'_i, z'_i) \ \& \ l(\delta) \geq x_{i+1}\}$ for $i \leq m$.
\end{itemize}

It remains to be checked that $x R_n^\al w$. Clearly, $x_i \geq w_i$ for $i \, {>} \, n$. Also, since $w_n = y_n$, $w_{n+1} = y_{n+1}$ and $x_n \geq y'_n = y_n + \big(1+ e(y_{n+1})\big) \cdot \al$ we have that $x_n \geq w_n + \big(1+ e(w_{n+1})\big) \cdot \al$. Thus, we need to see that $x_i \, {>} \, w_i$ for $i \, {<} \, n$.  For $i, \ m \, {<} \, i \, {<} \, n$, we can easily check that $y'_i \, {>} \, y_i = w_i$, and so $x_i \, {>} \, w_i$. For $i \leq m$, we show by induction on $k$ that $x_{m-k} \, {>} \, w_{m-k}$. For the base case, we can have that $x_{m+1} \, {>} \, w_{m+1}$. Also we can observe that $\max(y'_m, z'_m) \geq    \max(y_m, z'_m)$. Therefore $x_{m} \, {>} \, w_{m}$. For the inductive step, by the I.H. we have that $x_{m-k} \, {>} \, w_{m-k}$. Again, $\max(y'_{m - (k+1)}, z'_{m - (k+1)}) \geq \max(y_{m - (k+1)}, z'_{m - (k+1)})$, and so $x_{m-(k+1)} \, {>} \, w_{m-(k+1)}$. Hence, in virtue of Theorem \ref{theorem:CharacterizationOfRnAlphaSteps} we get that $x R_n^\al w$, that is, $x \Vdash \la n^\al \ra \big( \, \varphi \, \land \, \la m^{\be+ 1} \ra \psi \,  \big)$.
\end{proof}

Although it is not needed later in this paper, we find it useful to observe that for any $x = \la x_0, \ldots, x_k, 0 \ra \in I$ there is $\psi \in {\sf MNF}$ so that $\llbracket x \rrbracket = \llbracket \psi \rrbracket =  \llbracket x_\psi \rrbracket$. Having finite support is essential since e.g. the Ignatiev \mbox{sequence $\la \varepsilon_0, \varepsilon_0, \ldots \ra \in I$} is not modally definable.

\section{Completeness}

To establish the completeness of our system, first we need the following proposition that characterizes the non-derivability between formulas in monomial normal form.

\begin{proposition} \label{failComesFromPsi}
Given $\varphi, \, \psi \in {\sf MNF}$, if $\varphi \not \vdash \psi$ then there is $m_I \in {\sf N\text{-}mod}(\psi)$ such that $\pi_{m_I} (\varphi) < \pi_{m_I} (\psi)$.
\end{proposition}
\begin{proof}
Assume $\varphi \not \vdash \psi$ and suppose, towards a contradiction, that for any \linebreak $m \in {\sf N\text{-}mod}(\psi)$ we have that $\pi_m(\varphi) \geq \pi_m(\psi)$. Then, by Definition \ref{projection}, for all $m > \max {\sf N\text{-}mod}(\psi)$ we also have that $\pi_m(\varphi) \geq \pi_m(\psi)$. For $m \not \in  {\sf N\text{-}mod}$ with
$m < \max {\sf N\text{-}mod}(\psi)$, we can observe that $\pi_m(\psi) = e^k \big( \pi_{n_I} (\psi)\big)$ where $n_I$ is the least element in ${\sf N\text{-}mod}(\psi)$ such that $n_I > m$ and $k = n_I - m$. Since by supposition, $\pi_{n_I}(\psi) \leq \pi_{n_I}(\varphi)$ then $e^k \big(\pi_{n_I} (\psi)\big) \leq e^k \big( \pi_{n_I}(\varphi) \big)$, and so $\pi_m(\psi) \leq \pi_m(\varphi)$. Hence, we can conclude that for all $m < \omega$, we have that $\pi_m(\varphi) \geq \pi_m(\psi)$, and thus, by Theorem \ref{CharacterizatioMNFDeriv}, $\varphi \vdash \psi$ contradicting our assumption.
\end{proof}

\begin{corollary} \label{pointOfFailure}
For any $\varphi, \, \psi \in {\sf MNF}$, if $\varphi \not \vdash \psi$ then $x_\varphi \Vdash \varphi$ and $x_\varphi \not \Vdash \psi$.
\end{corollary}
\begin{proof}
By Lemma \ref{trueInTranslation}, Item (\ref{trueInTrans2}), $x_\varphi \Vdash \varphi$. On the other hand, by Proposition \ref{failComesFromPsi} and Lemma (\ref{trueInTranslation}), Item \ref{trueInTrans1} we have that $x_\varphi \not \Vdash \la m_I^{\be_I} \ra \top$ where $\be_I = \pi_{m_I}(\psi)$. Hence, $x_\varphi \not \Vdash \psi$. 
\end{proof}

With these tools, we can easily prove the completeness of $\TPr$.

\begin{theorem}[Completeness] \label{completeness}
Given formulas $\varphi, \, \psi \in \FLE$, if $\varphi \models \psi$, then $\varphi \vdash \psi$.
\end{theorem}
\begin{proof}
By Theorem \ref{MNFT}, w.l.o.g. let $\varphi, \, \psi \in {\sf MNF}$. Reasoning by contraposition, suppose $\varphi \not \vdash \psi$. Therefore, by Corollary  \ref{pointOfFailure}, $x_\varphi \Vdash \varphi$ but $x_\varphi \not \Vdash \psi$. Therefore,  $\varphi \not \models \psi$.
\end{proof}

\nocite{*}
\bibliographystyle{plain}
\bibliography{ref1}

\end{document}